\documentclass[10pt]{amsart}

\usepackage{amsmath,amssymb,amsthm,amsrefs}
\usepackage{latexsym}
\usepackage{indentfirst}
\usepackage{graphicx}
\usepackage{placeins}
\usepackage{booktabs}
\usepackage{algorithm}
\usepackage{algorithmic}
\usepackage{multirow}

\theoremstyle{plain}

\theoremstyle{definition}

\theoremstyle{remark}

\newcommand{\bbm}{\begin{bmatrix}}
\newcommand{\ebm}{\end{bmatrix}}

\newcommand{\R}{\mathbb{R}}

\newcommand{\T}{\mathsf{T}}

\renewcommand{\P}{\mathrm{P}}

\begin{document}

\title[Independent particle approximation to DPP]{A heuristic independent particle approximation to
  determinantal point processes}

\author[]{Lexing Ying}
\address[Lexing Ying]{Department of Mathematics and ICME, Stanford University, Stanford, CA 94305}
\email{lexing@stanford.edu}

\thanks{The work of L.Y. is partially supported by the U.S. Department of Energy, Office of Science,
  Office of Advanced Scientific Computing Research, Scientific Discovery through Advanced Computing
  (SciDAC) program and also by the National Science Foundation under award DMS-1818449.
}



\begin{abstract}
A determinantal point process is a stochastic point process that is commonly used to capture
negative correlations. It has become increasingly popular in machine learning in recent
years. Sampling a determinantal point process however remains a computationally intensive task. This
note introduces a heuristic independent particle approximation to determinantal point processes. The
approximation is based on the physical intuition of fermions and is implemented using standard
numerical linear algebra routines. Sampling from this independent particle approximation can be
performed at a negligible cost. Numerical results are provided to demonstrate the performance of the
proposed algorithm.
\end{abstract}

\maketitle

\section{Introduction}\label{sec:intro}

A determinantal point process is a stochastic point process that is commonly used to capture
negative correlations \cite{macchi1975coincidence}. Let $S$ be a set of points. A determinantal
point process is a random set $A$ such that the probability of $\{x_1,\ldots,x_k\}\subset A$ is
given by $\det (K(x_i,x_j))_{1\le i,j\le k}$, where $K: S\times S\rightarrow \R$ is symmetric
positive semi-definite operator. To simplify the discussion, we assume for simplicity that $S$ is a
discrete set of size $N$ and $A$ has a fixed size $k$. Such a DPP is called elementary and
\begin{equation} \label{eqn:dpp}
  \P(A = \{x_1,\ldots,x_k\}) = \det (K(x_i,x_j))_{1\le i,j\le k}.
\end{equation}
In what follows, we shall also refer to the diagonal $\rho(x)\equiv K(x,x)$ of $K$ as density.

Many natural point processes can be modeled by DPPs. Examples include distribution of
non-interacting Fermions, descent subsequences in random sequences, non-intersecting random walks,
edge distributions of spanning trees, and eigenvalue of random matrices
\cite{burton1993local,soshnikov2000determinantal,johansson2004determinantal,borodin2009determinantal,tao2009determinantal}. More
recently, DPP has played a significant role in improving fairness and diversity of sampling
algorithms in modern machine learning \cite{kulesza2012determinantal}.

However, sampling from DPP remains to be a challenging computational problem. In
\cite{hough2006determinantal} Hough et al proposed the standard DPP sampling algorithm with
$O(Nk^3)$ complexity.  Though various improvements (e.g. \cite{kulesza2010structured}) and
approximate algorithms (e.g. \cite{deshpande2010efficient}) have been proposed, sampling from DPP
remains a hard computation problem.

\subsection{Problem and contribution}
This note considers the problem of approximating DPP with an independent particle process with
disjoint support, i.e., generating the $k$ samples by sampling each from a support disjoint
region. More specifically, we introduce a disjoint union $S = S_1 \cup \ldots \cup S_k$ and
associate with each $S_i$ a probability density $\rho_i(\cdot)$ supported on $S_i$. Each realization
of this approximate process is then generated by sampling one point from each
$\rho_i(\cdot)$. Equivalently, this also amounts to approximating the kernel matrix $K$, after an
appropriate reordering, with a block-diagonal matrix, where each diagonal block has rank one and
unit trace. The main advantage of this approximation is that sampling from this approximation can be
done extremely rapidly. 

The main contribution of this note is to introduce a simple heuristic algorithm for constructing
such an approximation. The algorithm itself requires no more than standard numerical linear algebra
routines and the numerical results are provided to demonstrate its effectiveness.

\subsection{Motivation.}

It is natural to ask why one could expect such an approximation to be reasonable. This
approximation, though crude sometimes, is well-motivated from physics and chemistry
consideration. In the work of Macchi \cite{macchi1975coincidence}, DPP is originally named {\em
  fermionic point process} and the elementary DPP with $k$ points is exactly the distribution
function of $k$ non-interacting electrons. In this language, $K(x,x')$ is the density matrix and
$\rho(x)$ is the single-electron density. Here, {\em non-interacting} means that there is no further
interaction between the electrons besides the Pauli's exclusion principle
\cite{negele2018quantum}. To see this, note that the matrix $K$ can be decomposed as
\begin{equation}\label{eqn:Kphi}
  K(x,x') = \sum_{i=1}^k \phi_i(x) \phi_i(x')
\end{equation}
where the functions $\{\phi_i(x)\}$ from $S$ to $\R$ are called electron orbitals. The
multi-electron wave function $\Phi(x_1,\ldots,x_k)$ is given by the {\em Slater determinant}
\cite{slater1929theory,lin2019mathematical}
\[
\Phi(x_1,\ldots,x_k) = \frac{1}{\sqrt{k!}} \det(\phi_i(x_j))_{1\le i,j\le k}.
\]
The multi-electron density is then the square of the wave function,
\[
|\Phi(x_1,\ldots,x_k)|^2 = \frac{1}{k!} |\det(\phi_i(x_j))|^2 = \frac{1}{k!} \det(K(x_i,x_j))_{1\le
  i,j\le k}.
\]
Since electrons are indistinguishable, the probability of finding the $k$ electrons at the location
set $\{x_1,\ldots,x_k\}$ is given by
\[
P(\{x_1,\ldots,x_k\}) = k! \cdot \frac{1}{k!} \det(K(x_i,x_j)) = \det(K(x_i,x_j)),
\]
which matches \eqref{eqn:dpp} exactly.

The choice of the orbitals $\{\phi_i(x)\}$ in \eqref{eqn:Kphi} is not unique: applying an arbitrary
$k\times k$ orthogonal matrix to $\{\phi_i(x)\}$ generates an equally valid set of orbitals and
keeps $K$ and the DPP unchanged. However, different sets of orbitals do have different physical
interpretations and computational implications. It is often preferred to choose a set of localized
orbitals, such as the atomic orbitals $1s, 2s, 2p_x, 2p_y, 2p_z, \ldots$. With these localized
orbitals identified, one often says that a certain electron is in a certain orbital, without
actually referring to the actual position of the other electrons. This implies that, up to a
reasonable approximation, it makes sense to sample each electron from its own orbital independently
due to their locality. However, such a strategy can violate the Pauli exclusion principle, since the
localized orbitals can still overlap. As a result, in order to implement this intuition in our DPP
approximation, the algorithm needs to construct disjoint supports for each particle.

\subsection{Content.}
The rest of the note is organized as follows. Section 2 describes the algorithm for approximating
DPP with independent particles. Several numerical examples are presented in Section 3 to demonstrate
the performance of this heuristic algorithm. Finally, Section 4 includes some discussion for future
work.

\subsection{Data availability statement.}
Data sharing not applicable to this article as no datasets were generated or analyzed during the
current study.

\section{Algorithm} \label{sec:algo}

This section describes the algorithm for approximating DPP with independent particles. Let us recall
that, given a kernel matrix $K\in \R^{N\times N}$ defined on $S$, the objective is to introduce a
disjoint union $S = S_1 \cup \ldots \cup S_k$ and associate to each $S_i$ a probability density
$\rho_i(\cdot)$ supported within $S_i$. The algorithm consists of two steps: localization and
partitioning.

\subsection{Localization}
Let us introduce the matrix $\Phi=(\phi_1,\ldots,\phi_k) \in \R^{N\times k}$ with columns equal to
the orbitals $\{\phi_i(x)\}$ in \eqref{eqn:Kphi}. The task of this first step is to find an
equivalent set of orbitals $\{v_i(x)\}$, which are as localized as possible. Similar to the
definition of $\Phi$, we also introduce the matrix $V=(v_1,\ldots,v_k)\in\R^{N\times k}$.

Here, we follow the method of selected columns of density matrix (SCDM) introduced in
\cite{damle2015compressed}. In the matrix form, the task is to find an orthogonal matrix $O$ such
that
\[
V:=\Phi O
\]
has columns as localized as possible. Directly optimizing a locality measure/functional for $V$ over
all possible orthogonal matrices is a non-trivial optimization problem. Instead, the key idea of
SCDM is that the density matrix $K = \Phi \Phi^\T$ often has localized columns
\cite{benzi2013decay,benzi2007decay}. Therefore, instead of searching $O$ from the infinite set of
$k\times k$ orthogonal matrices, one can simply look for the columns of $O$ from the columns of
$\Phi^\T$. This can be implemented for example by performing a pivoted QR factorization to the
matrix $K=\Phi\Phi^\T$. However, since $\Phi$ is a matrix with orthogonal columns, it is equivalent
to perform the pivoted QR factorization to the smaller matrix $\Phi^\T$ instead:
\begin{equation} \label{eqn:qr}
  [Q,R,\sigma] = \text{qr}(\Phi^\T),
\end{equation}
where $Q$ is a $k\times k$ orthogonal matrix, $R$ is a $k\times N$ upper-triangular matrix, and
$\sigma$ is a $k$-dimensional integer vectors that identifies the first $k$ pivoted columns of
$\Phi^\T$. Given $\sigma$, the columns in $(\Phi \Phi^\T)(:,\sigma) = K(:,\sigma)$ are all localized
since they are selected columns of $K$. However, they are not orthogonal. To regain orthogonality,
one can set
\begin{equation} \label{eqn:V}
  V = K(:,\sigma) (K(\sigma,\sigma))^{-1/2}.
\end{equation}
It is easy to check that $V$ is indeed a matrix with orthonormal columns:
\begin{align*}
  V^\T V &= (K(\sigma,\sigma))^{-1/2} K(\sigma,:) K(:,\sigma) (K(\sigma,\sigma))^{-1/2} \\
  & = (K(\sigma,\sigma))^{-1/2} K(\sigma,\sigma) (K(\sigma,\sigma))^{-1/2} = I,
  \end{align*}
where the second step uses the fact that $K$ is a projection.

Though the definition of $V$ in \eqref{eqn:V} involves the $N\times N$ matrix $K$, it is equivalent
to write
\[
V = \Phi \Phi^\T(:,\sigma) (\Phi(\sigma,:)\Phi^\T(:,\sigma))^{-1/2},
\]
where it is clear that the computation of $V$ can be performed without any explicit reference to the
full matrix $K$.

\subsection{Partitioning}
The task of the second step is to partition $S$ into $k$ disjoint subsets $S_1,\ldots,S_k$. The
easiest way is to simply set
\begin{equation}\label{eqn:S}
  S_i = \{x\in S | i = \text{argmax}_j |v_j(x)|\},
\end{equation}
with random tie-breaking at a point $x\in S$ whenever multiple $v_j(\cdot)$ vectors have the same
absolute value at $x$.  Once $S_i$ is identified, one set the density $\rho_i(x)$ within each $S_i$
as
\begin{equation}\label{eqn:rho}
  \rho_i(x) =
  \begin{cases}
    \frac{\rho(x)}{\sum_{x\in S_i} \rho(x)}, & x\in S_i\\
    0, & x\not\in S_i.
  \end{cases}
\end{equation}
The main shortcoming of this approach is that the sum $\sum_{x\in S_i} \rho(x)$ can deviate
noticeably from $1$. Therefore, after the renormalization step in \eqref{eqn:rho}, the density
$\rho_i(x)$ can differ significantly from the original density $\rho(x)$.

To fix this issue, a heuristic balancing step is introduced. We seek for a set of scaling factors
$\{\alpha_i\}$ close to one such that the sets $S_i$ defined via
\begin{equation}\label{eqn:Sn}
  S_i = \{x\in S | i = \text{argmax}_j |(\alpha_jv_j)(x)| \}
\end{equation}
satisfy the constraints that for each $i$
\[
\sum_{x\in S_i} \rho(x) = 1.
\]
Once $\{\alpha_i\}$ are identified, one simply set
\begin{equation}\label{eqn:rhon}
  \rho_i(x) =
  \begin{cases}
    \rho(x), & x\in S_i\\
    0, & x\not\in S_i.
  \end{cases}
\end{equation}

\subsection{Sampling}
Once $\{S_i\}$ and $\{\rho_i\}$ are computed, sampling from this independent particle model is
straightforward.
\begin{enumerate}
\item For each $i=1,\ldots, k$, sample $x_i$ from $S_i$ following the distribution $\rho_i(x)$.
\item Return the set $\{x_1,\ldots,x_k\}$.
\end{enumerate}
$\{x_i\}$ are clearly disjoint since $\{S_i\}$ are disjoint. The cost of this sampling algorithm is
also extremely low. By adopting a binary search structure for the weights of $\rho_i(x)$, each
sample $x_i$ can be generated in $O(\log N)$ steps. Therefore, the overall sampling cost is $O(k
\log N)$.

\section{Numerical results} \label{sec:res}

This section considers several geometric sampling problems in two dimensional spaces to illustrate
the performance of the proposed heuristic algorithm. In each example, the functions $\{\phi_i(x)\}$
of $K(\cdot,\cdot)$ are given as input.

\begin{figure}[h!]
  \centering{
    \begin{tabular}{ccc}
    \includegraphics[width=0.32\textwidth]{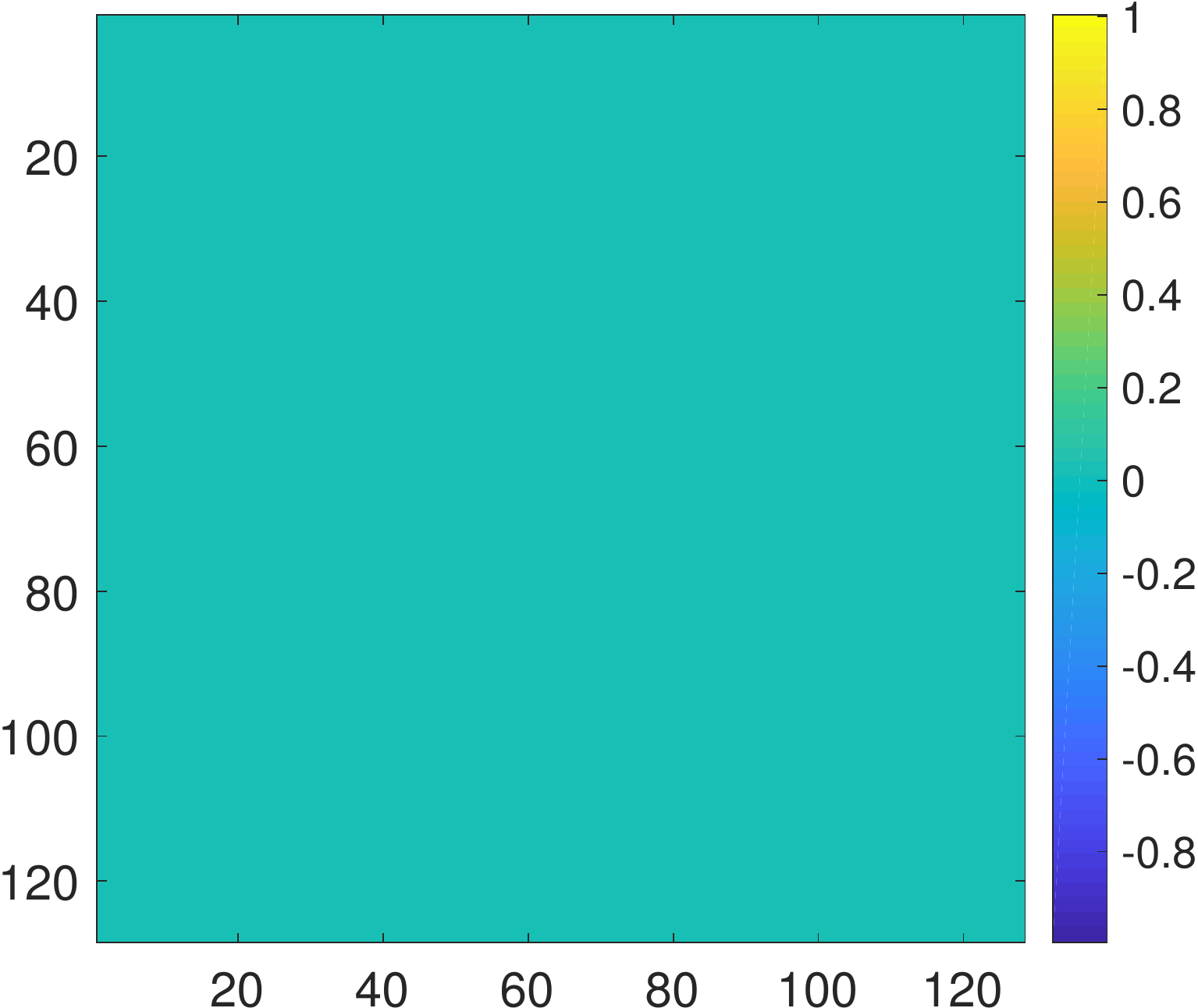} &
    \includegraphics[width=0.32\textwidth]{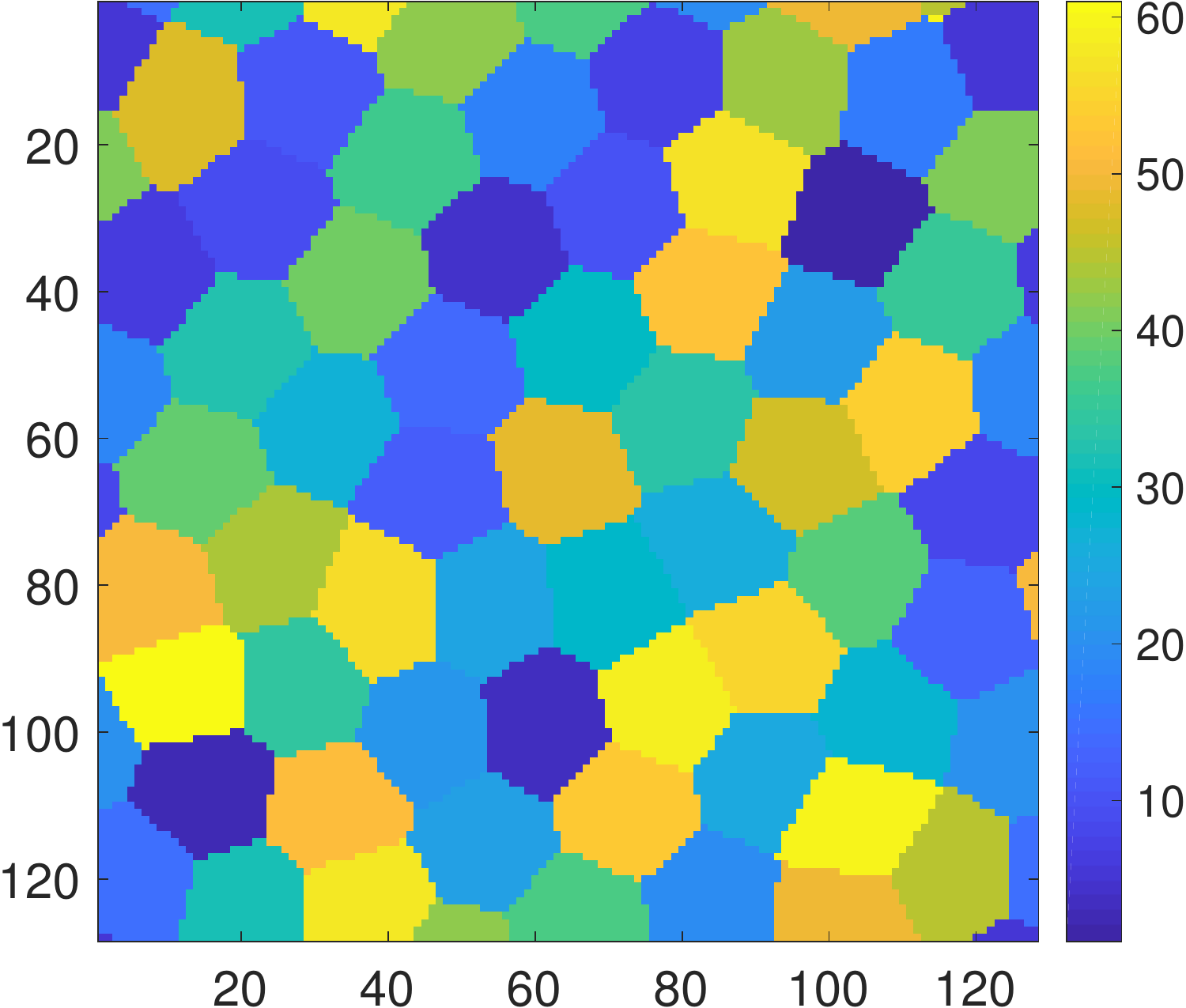} &
    \includegraphics[width=0.32\textwidth]{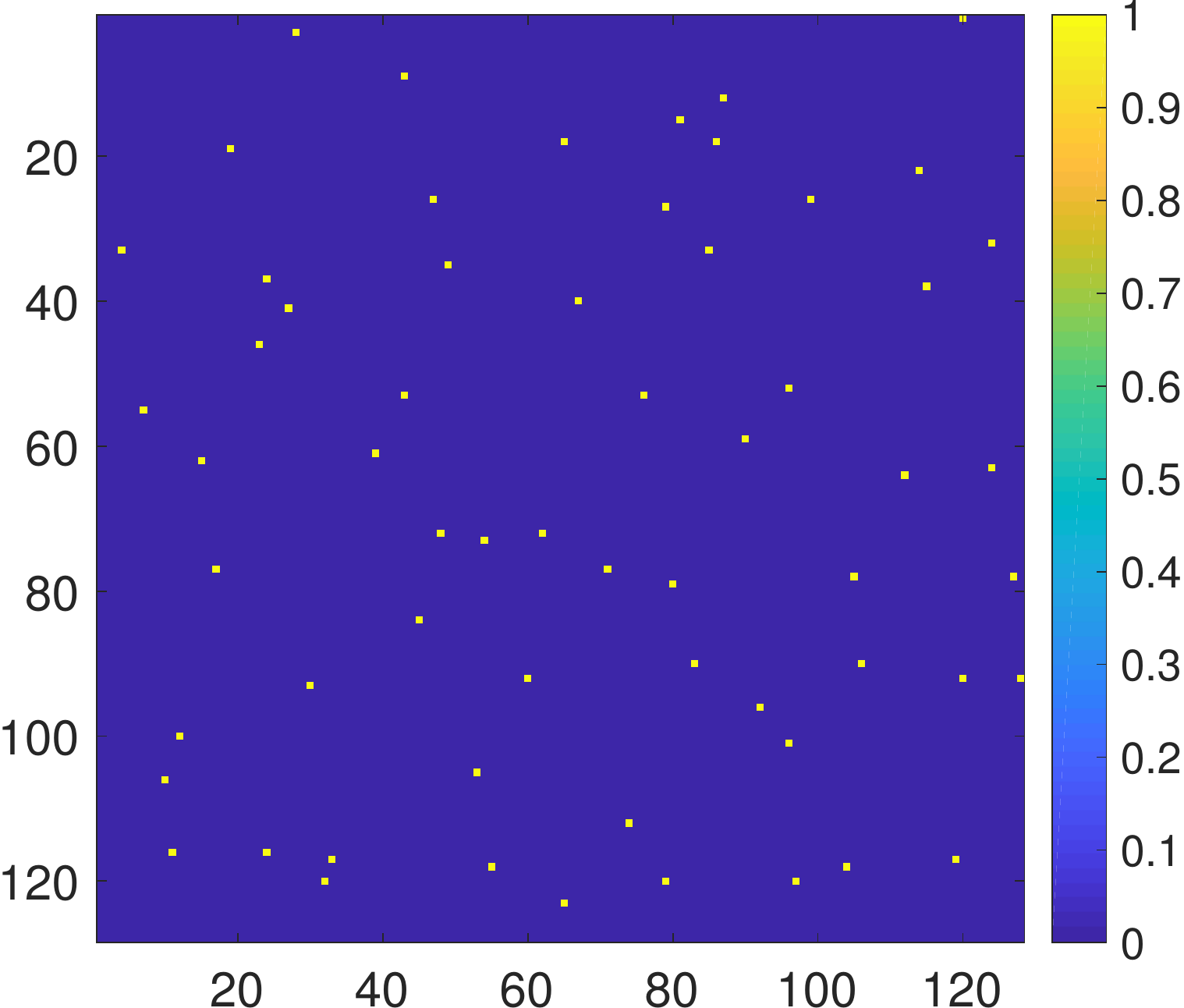}\\
    (a) & (b) & (c)
    \end{tabular}
  }
  \caption{The first example. (a) density $\rho(x)$. (b) partitioning $\{S_i\}$. (c) a realization
    of the resulting independent particle process.}
  \label{fig:uni}
\end{figure}

In the first example (see Figure \ref{fig:uni}), $S$ is a uniform Cartesian grid of $[0,1]^2$ with
$\sqrt{N}=128$ points in each dimension. Clearly, $N=128^2$. We set $k=61$ and the orbital functions
$\{\phi_i(x)\}$ to be the $k=61$ lowest eigenmodes of the discrete Laplacian $-\Delta$ on $S$ with
periodic boundary condition. The kernel $K(x,x')$ is given by \eqref{eqn:Kphi} and in this case the
density $\rho(x)$ is a constant function. Figure \ref{fig:uni}(a) shows the density $\rho(x)$ on the
Cartesian grid. Figure \ref{fig:uni}(b) plots the supports of different regions $\{S_i\}$ with
different colors. This plot demonstrate that $\{S_i\}$ are highly localized due to the locality of
the new orbitals $\{v_i(x)\}$. Finally, Figure \ref{fig:uni}(c) shows one realization of the
resulting independent particle process.

\begin{figure}[h!]
  \centering{
    \includegraphics[width=0.32\textwidth]{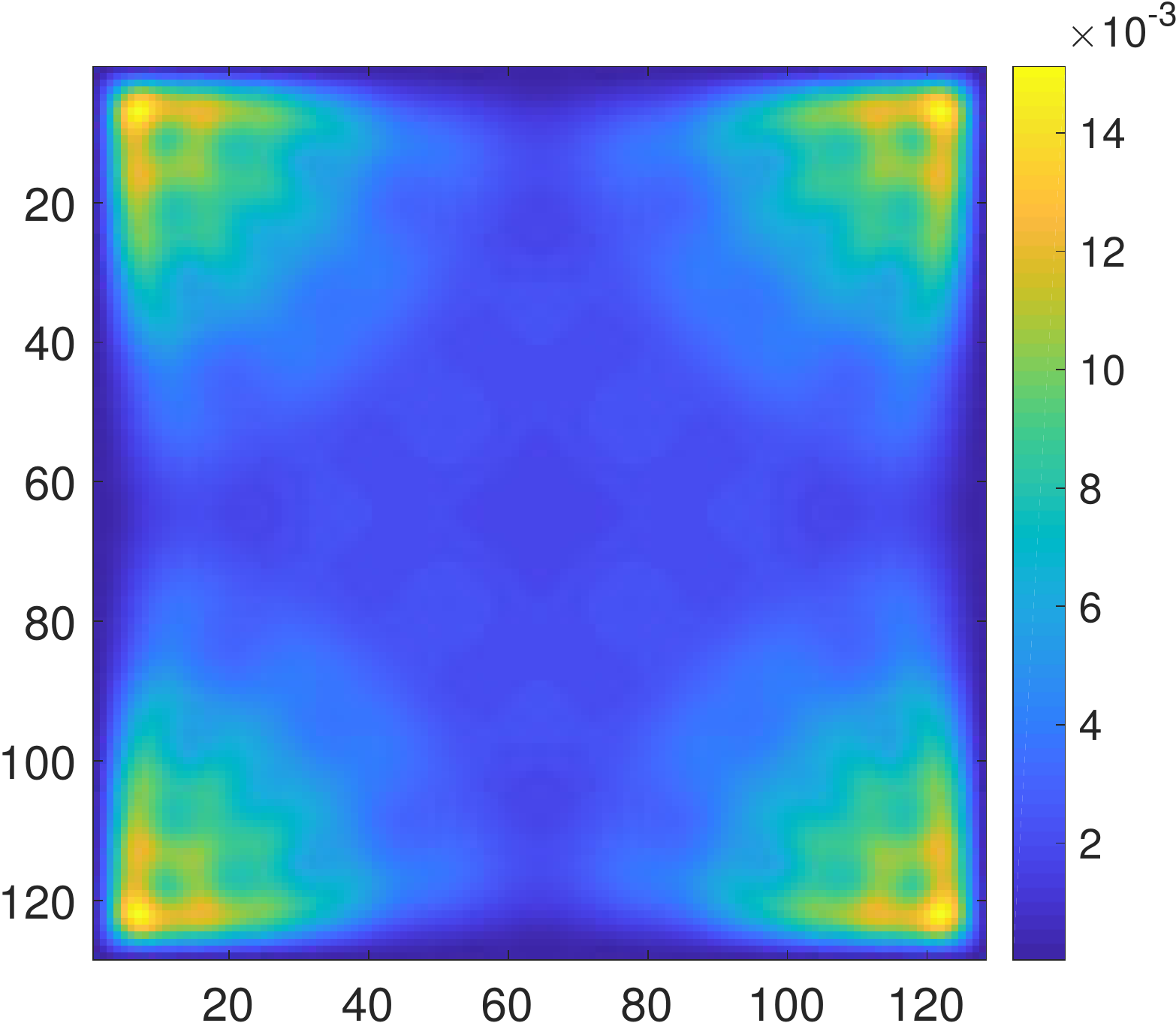} 
    \includegraphics[width=0.32\textwidth]{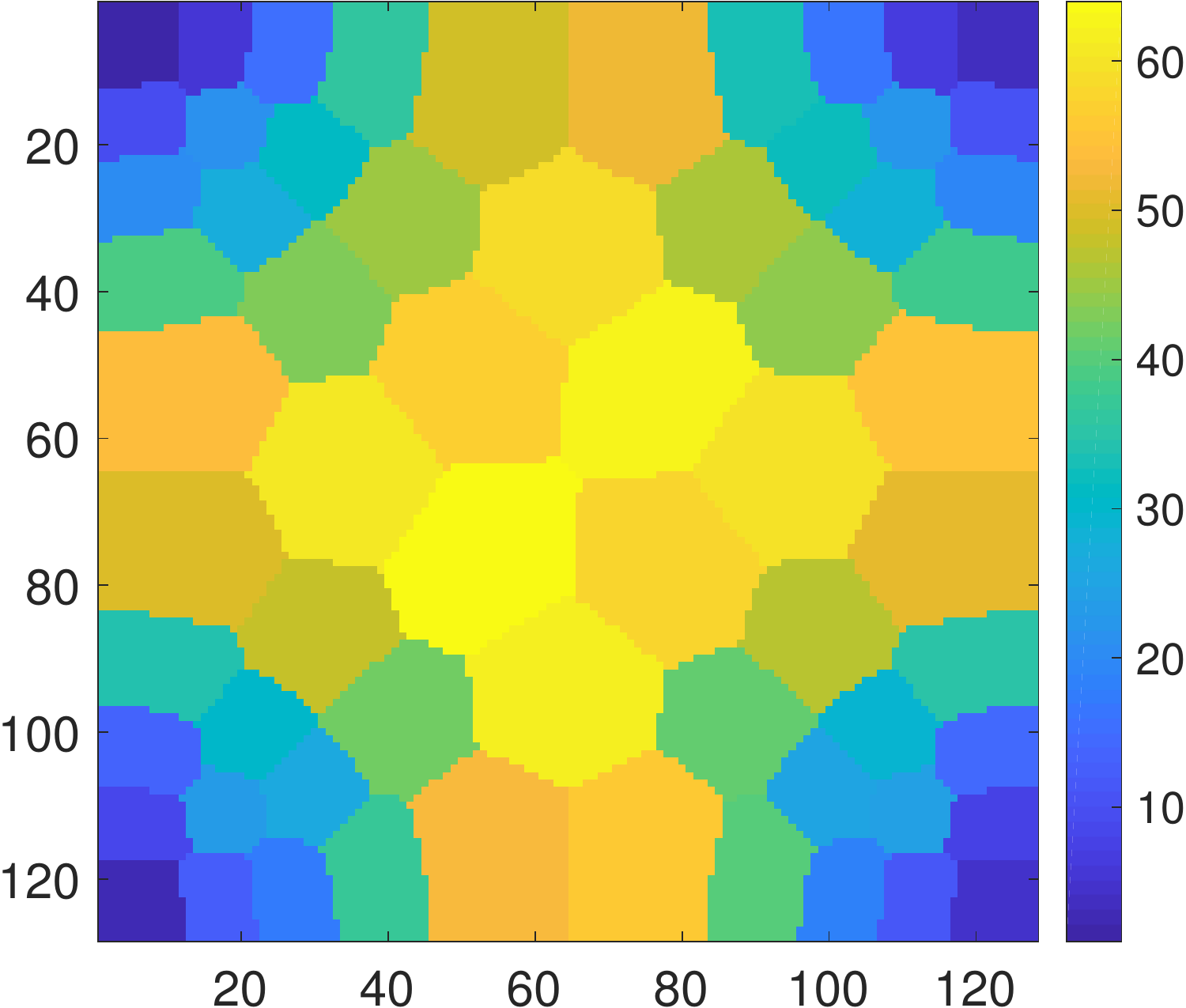} 
    \includegraphics[width=0.32\textwidth]{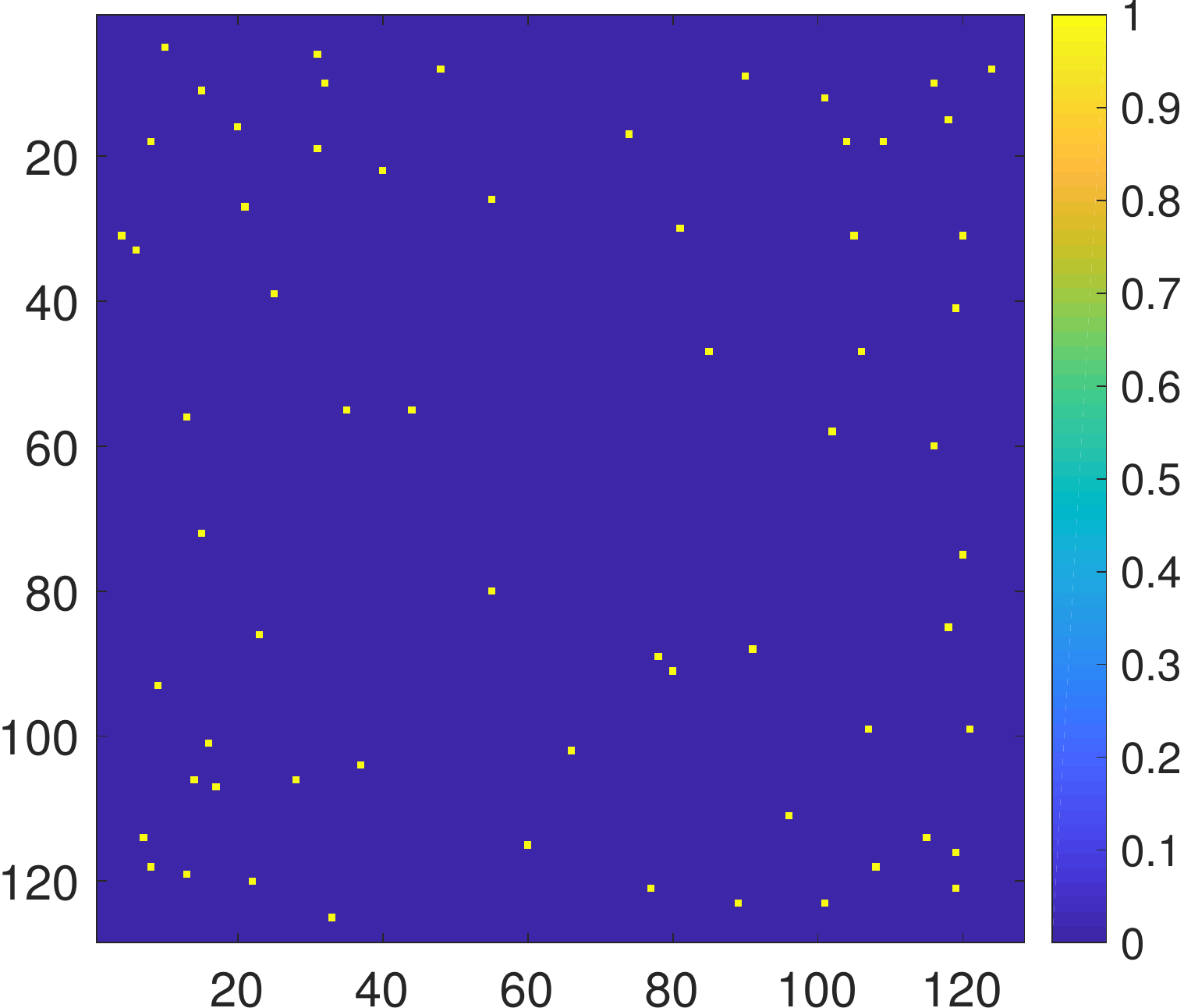}
  }
  \caption{The second example. (a) density $\rho(x)$. (b) partitioning $\{S_i\}$. (c) a realization
    of the resulting independent particle process.}
  \label{fig:ctr}
\end{figure}

In the second example (see Figure \ref{fig:ctr}), $S$ is the same uniform Cartesian grid and
$N=128^2$. We set $k=64$ and choose functions $\{\phi_i(x)\}$ to be the lowest eigenmodes of the
differential operator
\[
-\Delta + U(x),\quad U(x)\equiv U((x_1,x_2)) = -512\cdot (\cos(2\pi x_1)+1)\cdot (\cos(2\pi x_2)+1),
\]
with zero boundary condition. The kernel $K(x,x')$ can again be obtained from \eqref{eqn:Kphi}. In
this case the density $\rho(x)$ grows significantly nearly the four corners due to the low potential
values there (see Figure \ref{fig:ctr}(a)). Figure \ref{fig:ctr}(b) demonstrates the supports of
$\{S_i\}$ with different colors. We see that $\{S_i\}$ are again highly localized with necessary
area changes in order to accommodate the density variation across the domain. Finally, Figure
\ref{fig:ctr}(c) provides one realization from the resulting independent particle process.

\begin{figure}[h!]
  \centering{
    \includegraphics[width=0.32\textwidth]{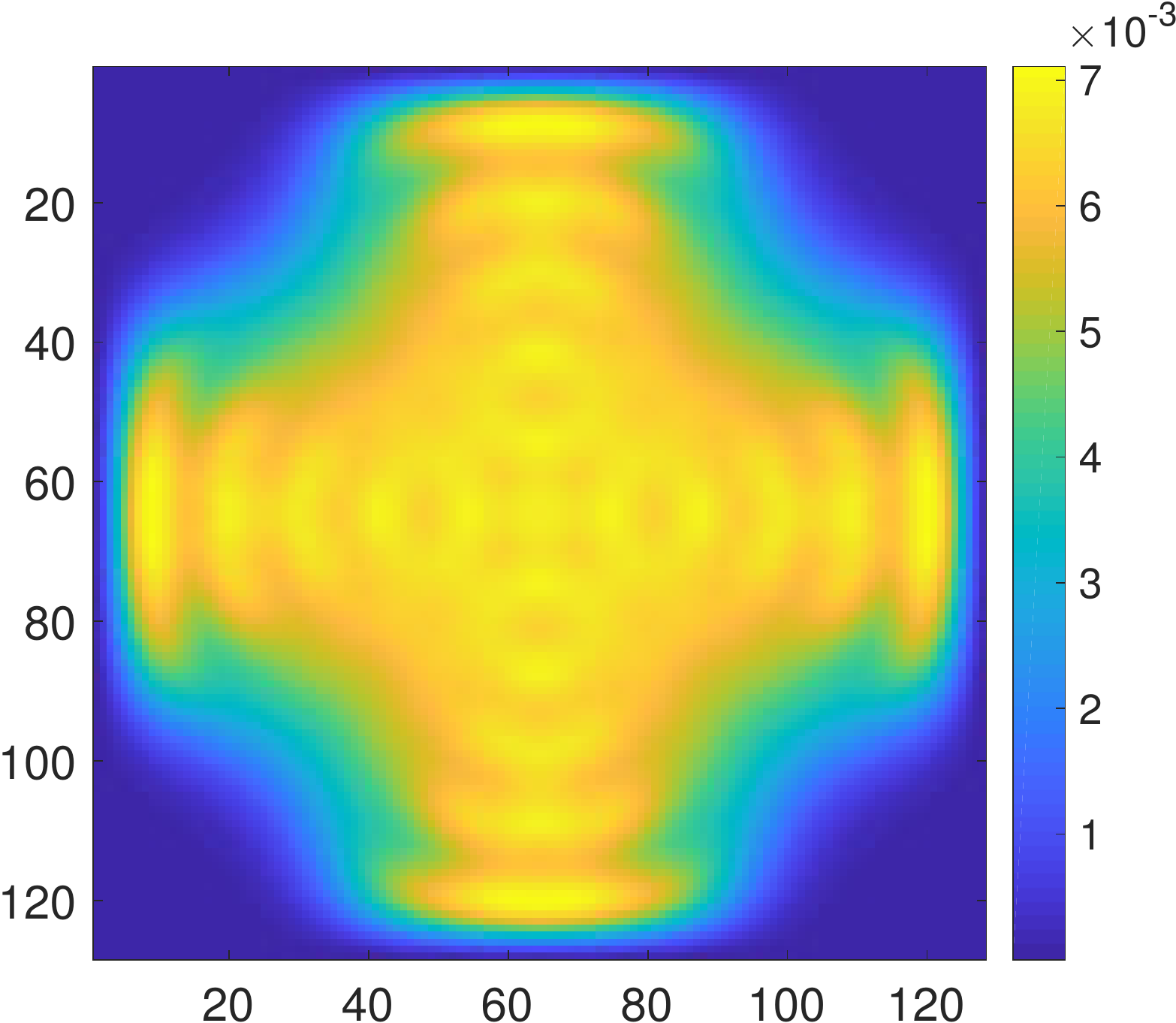} 
    \includegraphics[width=0.32\textwidth]{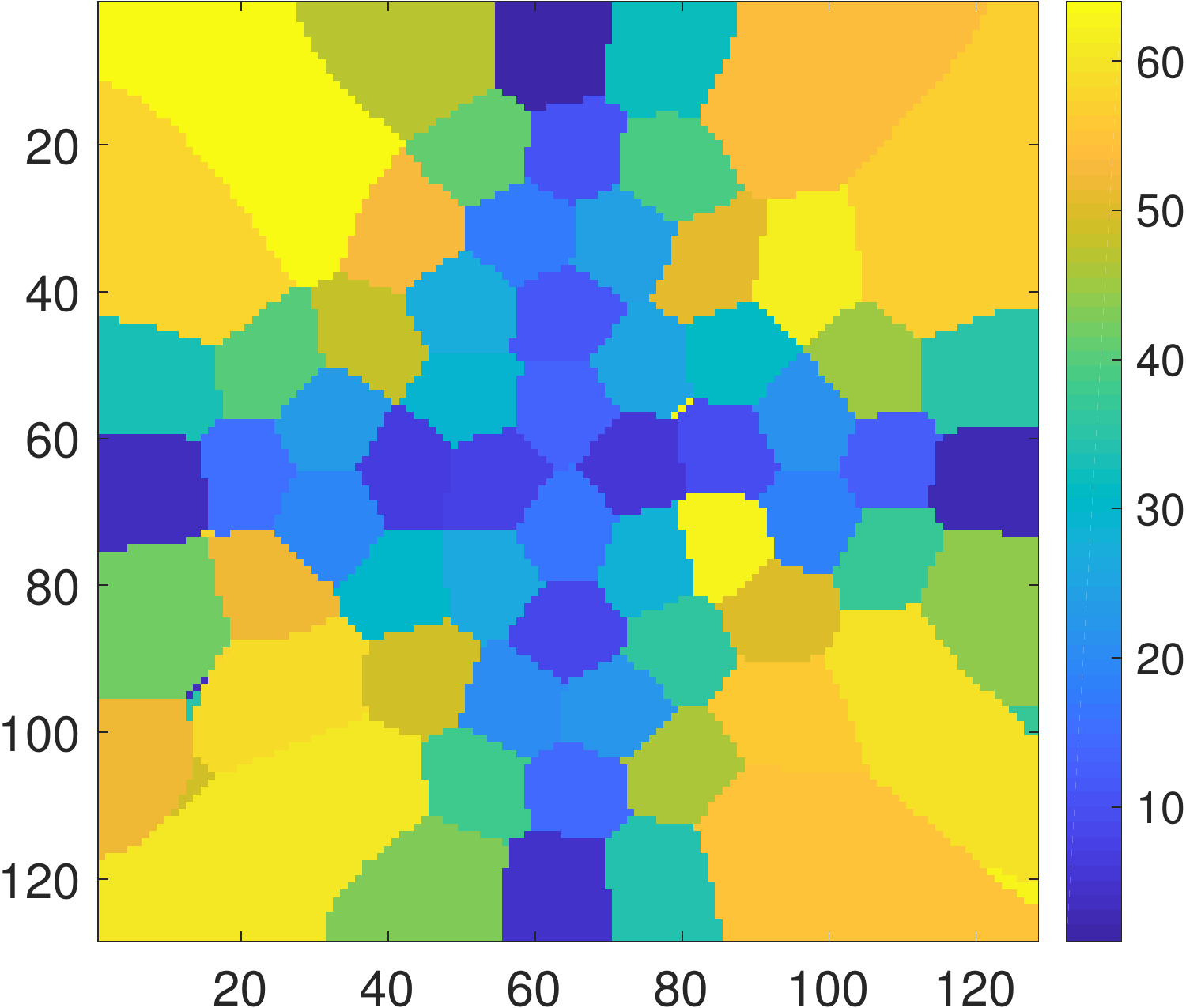} 
    \includegraphics[width=0.32\textwidth]{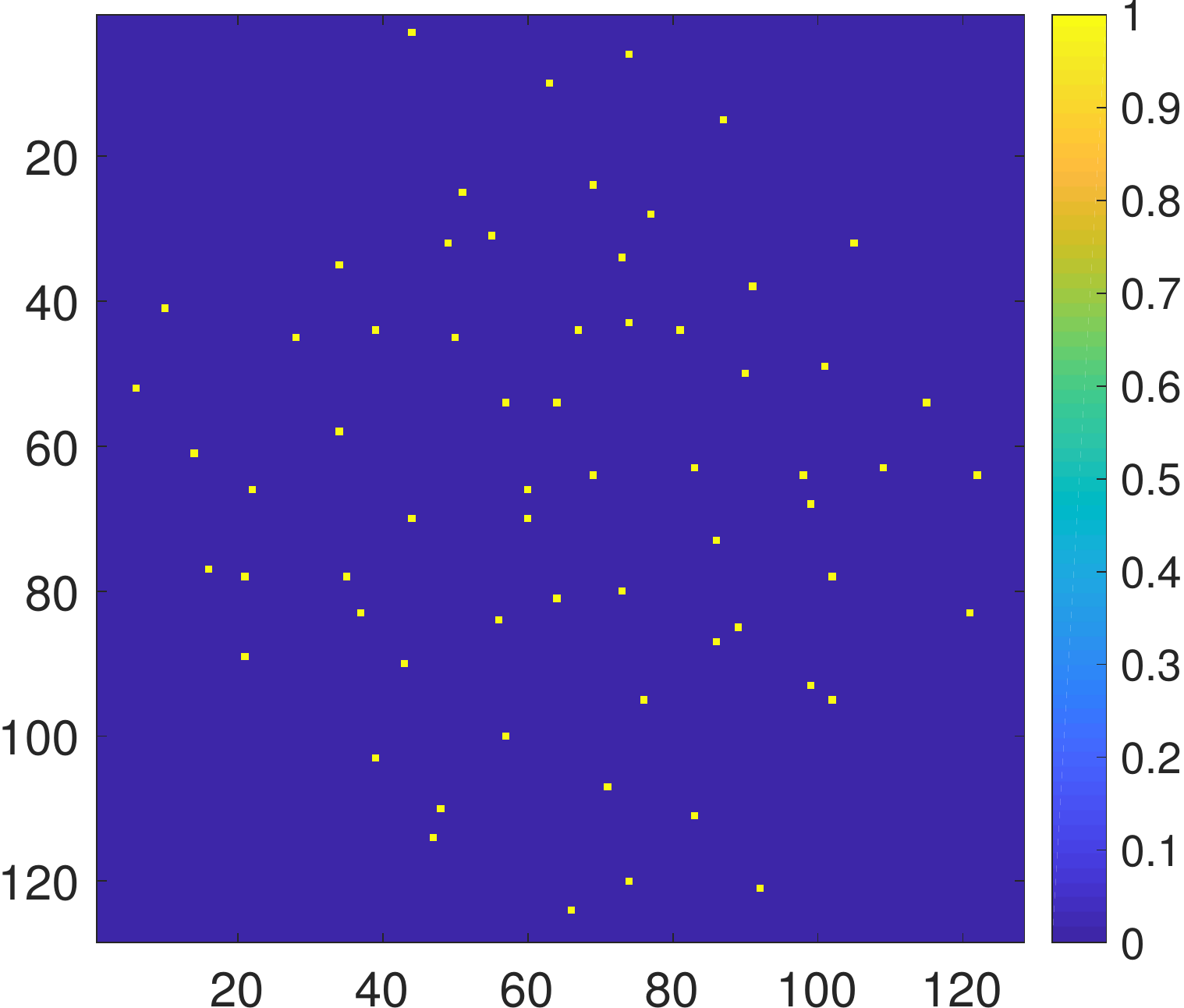}
  }
  \caption{The third example. (a) density $\rho(x)$. (b) partitioning $\{S_i\}$. (c) a realization
    of the resulting independent particle process.}
  \label{fig:cnr}
\end{figure}

In the third example (see Figure \ref{fig:cnr}), $S$ is still the uniform Cartesian grid with
$N=128^2$. We choose $k=64$ and let $\{\phi_i(x)\}$ be the lowest eigenmodes of the differential
operator
\[
-\Delta + U(x),\quad U(x)\equiv U((x_1,x_2)) = +512\cdot (\cos(2\pi x_1)+1)\cdot (\cos(2\pi x_2)+1),
\]
with zero boundary condition. The density $\rho(x)$ grows at the domain center due to the low
potential $U(x)$ there (see Figure \ref{fig:cnr}(a)). Figure \ref{fig:cnr}(b) demonstrates the
supports of $\{S_i\}$ with different colors.  $\{S_i\}$ are again highly localized with necessary
area changes to accommodate the density variation. Finally, Figure \ref{fig:cnr}(c) gives one
realization of the point set from this independent particle process.

\section{Discussions} \label{sec:disc}
This note introduces a heuristic independent particle approximation to determinantal point
processes.  The main benefit of this approximation is that it can be sampled with negligible cost.
There are several immediate directions for future work. First, this note only considers the
elementary DPP case, and it will be important to generalize this to general DPPs. Second, it will be
useful to explore the applications of this algorithm in machine learning applications where the
sampling speed of DPP is essential.

\bibliographystyle{abbrv}

\bibliography{ref}

\end{document}